\documentclass[12pt,reqno]{amsart}
\usepackage{amssymb}
\makeatletter
\def\proof{\@ifnextchar[{\@oproof}{\@nproof}}
\def\@oproof[#1][#2]{\trivlist\item[\hskip\labelsep
\textit{#2 Proof of\ #1.}~]\ignorespaces}
\def\@nproof{\trivlist\item[\hskip\labelsep\textit{Proof.}~]\ignorespaces}

\makeatother

\setbox0=\hbox{$+$}
\newdimen\plusheight
\plusheight=\ht0
\def\+{\;\lower\plusheight\hbox{$+$}\;}

\setbox0=\hbox{$-$}
\newdimen\minusheight
\minusheight=\ht0
\def\-{\;\lower\minusheight\hbox{$-$}\;}

\setbox0=\hbox{$\cdots$}
\newdimen\cdotsheight
\cdotsheight=\plusheight
\def\cds{\lower\cdotsheight\hbox{$\cdots$}}

\setbox0=\hbox{$+$}

\renewcommand{\(}{\left\(}
\renewcommand{\)}{\right\)}
\renewcommand{\[}{\left[}

\numberwithin{equation}{section}
 \theoremstyle{plain}
\newtheorem{theorem}{Theorem}[section]

\begin{document}
\title[A Study of Kummer's Proof of Fermat's Last Theorem for Regular Primes] {A Study of Kummer's Proof of Fermat's Last Theorem for Regular Primes}
\author{Manjil P.~Saikia}
\address{Department of Mathematical Sciences, Tezpur University, Napaam, Sonitpur, Pin-784028, India}
\email{manjil\_msi09@agnee.tezu.ernet.in, manjil.saikia@gmail.com}


\vspace*{0.5in}
\begin{center}
{\bf A Study of Kummer's Proof of Fermat's Last Theorem for Regular Primes}\\[5mm]
{\footnotesize  MANJIL P.~SAIKIA\footnote{Summer Research Fellow, Indian Institute of Science Education and Research (IISER), Sector-81, Mohali-140055, India.}}\\[3mm]

MATS137\\[3mm]

\emph{Summer Project under Prof.~Kapil Hari Paranjape.}
\end{center}
\vskip 5mm \noindent{\footnotesize{\bf Abstract.} We study Kummer's approach towards proving the Fermat's last Theorem for regular primes. Some basic algebraic prerequisites are also discussed in this report, and also a brief history of the problem is mentioned. We review among other things the Class number formula, and use this formula to conclude our study.}

\vskip 3mm


\vskip 3mm

\noindent {\footnotesize 2000 Mathematical Reviews Classification
Numbers: 11M06, 11M38, 11R18, 11R29, 11R37, 11R42.}

\section{{Introduction}}

One of the famous problems in all of mathematics is the Fermat's Last Theorem. A standard restatement is,

\begin{theorem}[Fermat's Last Theorem] There are no solutions to the following problem with $(X,Y,Z)$ integers

$$X^p+Y^p+Z^p=0$$

where $XYZ\neq0$ and $p\geq3$ is a prime.

\end{theorem}

The approach to the proof of Fermat's Last Theorem by Andrew Wiles can be thought of as a particular case of the following technique.

Suppose $(X,Y,Z)$ is a counter-example to Fermat's Last Theorem, then

\begin{itemize}

\item[(1)]To such a counter-example we attach a representation $$\rho(X,Y,Z):Gal(\overline{Q}/Q)\rightarrow GL_n(F_p)$$ Moreover, we have good ramification properties for this representation. For example,
    \begin{itemize}
    \item[(a)] The representation is unramified outside $p$,
    \item[(b)] The representation has "good" ramification properties at $p$.
    \end{itemize}

\item[(2)] The next step is to use Algebraic Number Theory to prove that such ramifications are not possible.
\end{itemize}

The proof of Kummer for the case of regular primes can also be thought in this light. In this report, we study the proof of Kummer of Fermat's Last Theorem for the case of regular primes. The approach of Kummer is to associates with every counter-example of Fermat's Last Theorem, a representation $$\rho:Gal(\overline{K}/K)\rightarrow F_p$$ where $K$ is a cyclotomic field of p-th roots of unity. Then, he gives us a way to finding out which primes $p$ are such that we have such a representation. Kummer showed that there are indeed such primes and his proof works only for 'regular' primes.

For the sake of completeness we give a brief explanation of regular primes. They may be defined in terms of the class number of cyclotomic fields or by means of Bernoulli numbers.

The Bernoulli numbers $B_0, B_1, B_2, \ldots$ are defined recursively: $$B_0=1$$ and for $n\geq 1$, $$\binom{n+1}{1}B_n+\binom{n+2}{2}B_{n-1}+\cdots+\binom{n+1}{n}B_1+1=0.$$

Thus, $B_1=\frac{-1}{2}$, $B_2=\frac{1}{6}$, $B_3=0$ and so on. It is easily seen that $B_{2k+1}=0$ for all $k\geq1$. The prime number $p$ is regular if $p$ does not divide the numerators of the Bernoulli numbers $B_2, B_4, \ldots, B_{r-5}, B_{r-3}$.

Let $p$ be an odd prime, let $$\zeta_p=\cos(2\pi/p)+i\sin(2\pi/p)$$ be a primitive root p-th root of $1$. Let $Q(\zeta_p)$ be the p-th cyclotomic field; it consists of all complex numbers of the form $$r_0+r_1\zeta_p+\cdots+r_{p-2}\zeta_p^{p-2},$$ with $r_o, r_1, \ldots, r_{p-1}\in Q$. The class number $h_p$ of $Q(\zeta_p)$ is a certain positive integer attached to $Q(\zeta_p)$; it is the number of classes of ideals of $Q(\zeta_p)$. Kummer showed that the prime $p$ is regular iff $p$ does not divide $h_p$. We shall view this later towards the end of the report.The smallest irregular prime is $37$. It is known that there are infinitely many irregular primes. On the other hand, it is conjectured, but it was never proved, that there are infinitely many regular primes.

 While working on the project, we reviewed many concepts and theorems of Algebra from \cite{dummit} and \cite{lang}. A few of those are:

 \begin{itemize}
 \item Quotient Field
 \item Cyclotomic Field
 \item Eisenstein Criterion
 \item Galois Group
 \item Galois Extension
 \item Ideals
 \item Dedekind Domains
 \item Class Group
 \item Ramification Properties
\end{itemize}

In section 2 we gave a very brief history of Fermat's Last Theorem till the 1847 proof by Kummer we are studying. Following this, in section 3, we recall some standard computations in the cyclotomic field of p-th roots of unity. In section 4, we show how a counter-example of Fermat's Last Theorem (if it exists) can be used to construct a cyclic extension of order $p$ of the cyclotomic field which is unramified every where. In section 5, we review the Class number formula, and then finally in section 6 we use this formula to check where such unramified extensions do indeed exist, concluding our study.

We fix a prime $p\geq5$ throughout the report, unless otherwise stated.

\section{{Fermat's Last Theorem: A Brief History}}
Sometime in the middle of the $17^{th}$ century Pierre de Fermat, an amateur French mathematician wrote in the margin of his copy of Bachet's edition of the works of Diophantus,
\begin{quote}It is impossible to separate a cube into two cubes, or a biquadrate
into two biquadrates, or in general any power higher than the second
into powers of like degree; I have discovered a truly marvelous proof,
which this margin is too small to contain.\end{quote}
This is the celebrated Fermat's Last Theorem, which in modern language translates to,
\begin{theorem}
If $n$ is any natural number greater than $2$, the equation $$X^n+Y^n=Z^n$$ has no solutions in integers, all different from $0$.\end{theorem}

Whatever marvelous proof Fermat had for his theorem nobody found out because in all of Fermat's letters to other mathematicians he never mentioned it. His theorem gathered much publicity after his death and remained the most prized unsolved problem of mathematics for more than three centuries. Euler, the legendary Swiss mathematician proved the theorem for $n=3$, and this was followed by Sophie Germain with her proof of the theorem for relatively small primes. The case when $n=5$ were proved with the efforts of Dirichlet and Legendre. Dirichlet further proved the theorem for $n=14$. Sometime before Kummer, the mathematician Lam\'{e} disposed off with the case $n=7$. This was then followed by Kummer's marvellous achievement which we shall review in the following few sections.

\section{{Arithmetic of Prime Cyclotomic Fields}}
Let $R$ denote the subring of complex numbers generated by $\omega=exp(2\pi i/p)$; let $K$ denote the quotient field of $R$, which is called the cyclotomic field of p-th roots of unity. In this section we review and state a few results that will be used in the next section about the ring $R$ and the field $K$.

The ring $R$ is isomorphic to $Z[X]/\Phi (X)$, where $$\Phi_p(X)=X^{p-1}+\ldots+X+1=\frac{X^p-1}{X-1}$$ is an irreducible polynomial. The irreducibility can be easily proved via the Eisenstein criterion. The field $K$ here is a Galois extension of $Q$ with Galois group $F_p^{\ast}$, which is a cyclic group of order $(p-1)$. We shall use $\gamma$ for a choice of generator, and use $\overline{\alpha}$ to denote $\gamma^{(p-1)/2}(\alpha)$ since $\gamma^{(p-1)/2}$ is the restriction of complex conjugation to $R$.

The ring $R$ is a Dedekind domain, whose prime ideals are
\begin{itemize}
\item If $q \in Z$ is a prime different from $p$, then let $f$ be the order of $q$ in $F_p^{\ast}$ and let $g=\frac{p-1}{f}$. Then there are $g$ prime ideals $Q_1, Q_2, \cdots, Q_g$ in $R$ such that their norms are $q^f$.

\item The element $\lambda = 1-\omega$ is prime in $R$ and $\lambda^{p-1}=(unit). p$.

\end{itemize}

We do not know a closed form expression of the generators of the group $U$ of units of $R$, but the numbers $$u_j=\gamma^{j}(\lambda)/\lambda=1+\omega+\cdots+\omega^{j-1}$$ are obviously in $R$ and are units there. The subgroup $U_{cycl}$ of the group $U$ of units of $R$ generated by $u_j$, $j=2, \ldots,(p-1)$ is called the group of cyclotomic units. If $u$ is a unit in $R$, then $\overline{u}/u$ is a root of unity in $R$. The roots of unity in $R$ are all of the form $\pm\omega^{j}$ for some $j=0,\ldots, p-1$. An element of $R$ is a p-th power only if it is congruent to an integer modulo $pR$. Hence, it follows that $\overline{u}/u=\omega^j$ for some $j$.

Let $L$ denote the subfield of $K$ fixed by complex conjugation; let $S=L\cap R$, then $L$ is a Galois extension of $Q$ with Galois group $F_p^{\ast}/\{\pm 1\}$. We see that $K$ is purely imaginary as no complex embeddings of $K$ have image within real numbers; while $L$ is real as all complex embeddings of $L$ have image within real numbers.

$S$ is a Dedekind domain and its ideals are described by:
\begin{itemize}
\item If $q \in Z$ is a prime different from $p$, then let $f\prime$ be the order of $q$ in $F_p^{\ast}/\{\pm 1\}$ and let $g\prime=\frac{p-1}{2f\prime}$. Then there are $g\prime$ prime ideals $Q_1, Q_2, \cdots, Q_{g\prime}$ in $R$ such that their norms are $q^{f\prime}$.

    \item The element $\mu=1-(\omega + \omega^{-1})$ is prime in $R$ and $\mu^{\frac{p-1}{2}}=(unit).p$.
    \end{itemize}

We have already seen that for any unit, $u$ in $R$ $\overline{u}/u=\omega^r$ for some integer $r$. Also $r\equiv 2s~(mod~p)$ for some integer $s$, so $u_1=\omega^{-s}u$ is in $S$. Thus, any unit in $R$ is the product of a root of unity and a unit in $S$.

 If $I$ is any ideal in $S$ then $IR$ is principal in $R$ iff $I$ is principal in $S$. So, the homomorphism from the class group of $S$ to that of $R$ is injective. In particular the order $h$ of the class group of $R$ is divisible by the order $h_{+}$ of the class group of $S$.

 If we have a unit $u$ in $R$,  such that it is congruent to an integer modulo $pR$ and if $u$ is itself not a p-th power, then the field extension of $K$ obtained by adjoining a p-th root of $u$ is a cyclic extension of $K$ of order $p$ which is unramified everywhere.

Again we have an important result from Class field theory that if there is an ideal $I$ in $R$ such that $I^p$ is principal and $I$ is not principal, then there is a cyclic extension of $K$ of order $p$ which is unramified everywhere. This follows from the identification of the class group of $R$ with the Galois group of the maximal unramified abelian extension of $K$. Then we use the fact that if an abelian group has an element of order $p$, then it has non-trivial character or order $p$.

\section{{Construction of Cyclic Order}}

Our aim is to show that if we have a counter-example to Fermat's Last Theorem, then there is a cyclic extension of order $p$ of $K$ which is unramified everywhere. We assume that the given counter-example $(X,Y,Z)$ has the property that these are mutually co-prime integrs.

We work out the details by dividing it into two cases. We shall be using most of the results from Section 4.

\textbf{Case 1: $p \nmid XYZ$}

First we see that $(X,Y,Z)$ are not all congruent modulo $p$. If not, we have $$3X\equiv X+Y+Z\equiv X^p+Y^p+Z^p\equiv 0~(mod~p).$$ Now, we assume that $p\geq 5$ and we get $X\equiv 0~(mod~p)$ which contradicts our hypothesis for this case.

Also, we see that $(X+\omega^j Y)$ are mutually co-prime in $R$ for $j=o,\ldots, p-1$. If not, then we have a prime ideal $P$ in $R$ containing $(X+\omega^j Y, X+\omega^k Y)$. Then this ideal $P$ contains $(1-\omega^{j-k})Y$. Now from the factorization $$(-Z)^p=X^p+Y^p=(X+Y)(X+\omega Y)\cdots(X+\omega^{p-1}Y)$$ we see that $P$ contains $Z$. Hence, by the assumption that $(X,Y,Z)$ are mutually co-prime we see that $P$ contains $(1-\omega^l)$ for some $0\leq 1\leq p-1$. By the previous section we see that $P=\lambda R$, but then $Z$ is a multiple of $p$ which contradicts our hypothesis in the present case.

From the above and unique factorization of ideals we see that we have ideals $I_j$ of $R$ such that $I_j^p=(X+\omega^jY)R$. We assume that $I_1$ is principal, then we have $$(X+\omega Y)=u.\alpha^p$$ for some $\alpha \in R$ and for some $u$, a unit in $R$. Applying the complex conjugation we obtain $$(X+\omega^{-1}Y)=\overline{u}.\overline{\alpha}^p$$.

By the results from the previous section we have $\omega^r\overline{u}=u$ for some $r$. Moreover, $\alpha^p$ is congruent to an integer modulo $pR$ and hence is congruent to its own complex conjugate. Thus we obtain $$X+\omega Y-\omega^rX-\omega^{r-1}Y\equiv 0~(mod~p).$$

From the results about $R$ in the previous section we see that it is a free abelian group with basis consisting of any $(p-1)$ elements of the set $\{1,\omega,\ldots,\omega^{p-1}\}$. From this and the fact that $X$ and $Y$ are prime to $p$ it follows easily that $r=1$ and $X\equiv Y~(mod~p)$.

By a similar reasoning and interchanging the roles of $Y$ and $Z$ we can conclude that there is an ideal $J_1$ such that $J_1^p=(X+\omega Z)$. Assuming $J_1$ is principal we see by an argument like the above that $X\equiv Z~(mod~p)$. However, these two congruences contradict the hypothesis in this case.

Hence, either $I_1$ or $J_1$ must be non-principal. But then from the result from Class field theory mentioned in the previous section we have the required cyclic extension of $K$.

\textbf{Case 2: $p\mid XYZ$}

We may assume that $Z=p^kZ_0$ and $(p,X,Y,Z_0$ are mutually co-prime, By writing $p=(unit).\lambda^(p-1)$ in the ring $R$, we obtain $$U^p+V^p+(unit)\lambda^{mp}W^p=0,~m>0$$ where $(U,V,W)$ are in $R$ so that $(U,V,W,\lambda)$ are mutually co-prime. Let $(U,V,W)$ be a collection of elements in $R$ that satisfy such an equation with $m$ the least possible. Then $\lambda$ divides one of the factors $(U+\omega^jV)$. But then we have $$(U+\omega^jV)-(U+\omega^kV)=\omega^j(1-\omega^{k-j})V=(unit).\lambda V$$ and thus, $\lambda$ divides all the factors $(U+\omega^jV)$. Moreover, since $V$ is co-prime to $p$ and thus $\lambda$ as well, we see that $(U+\omega^jV)/\lambda$ have distinct residue classes modulo $\lambda R$. But then, by the pigeon-hole principle there is at least one $0\leq j \leq (p-1)$ such that $(U+\omega^jV)$ is divisible by $\lambda^2$ in $R$. Replacing $V$ by $\omega^jV$ we may assume that $(U+V)$ is divisible by $\lambda^l$ for some $l>1$. Hence we may write $$U+V=\lambda^la_0$$ and $$U+\omega^kV=\lambda a_k; k>0$$ where all the $a_k$ are elements of $R$ that are co-prime to $\lambda$ and with each other like the previous case. This gives us the identity $l+(p-1)=mp$ or equivalently $l=(m-1)p+1$. Since, $l\geq 2$, we have $m\geq 2$.

Now by the unique factorization of ideals in $R$ we see that there are ideals $I_j$ in $R$ such that $I_j^p=a_jR$. We assume that $I_0, I_1$ and $I_{p-1}$ are principle, then we have the following $$U+V=\lambda^l.u.b_0^p$$ $$U+\omega V=\lambda.v.b_1^p$$ and $$U+\omega^{-1}V=\lambda.w.b_{-1}^p$$ for some units $u,v$ and $w$ in $R$ and some elements $b_0,b_1$ and $b_{-1}$ in $R$. Eliminating $U$ and $V$ from the above we get $$\lambda^l.u.b_0^p-\lambda.v.b_1^p=\omega(\lambda.w.b_{-1}^p-\lambda^l.u.b_0^p)$$ which becomes $$b_1^p+v_1.b_{-1}^p+\lambda^{l-1}.v_2.b_0^p=0$$ where $v_1$ and $v_2$ are units which can be justified by the fact that $1+\omega$ is a unit in $R$. Modulo $pR$ the last term on the left-hand side vanishes since $l\geq p>(p-1)$. Thus, we see that $v_1$ is congruent to a p-th power and thus an integer modulo $pR$. By the previous section we have a representation of Galois as required, unless $v_1$ is a p-th power. If $v_1=v_3^p$, then $U,V,W)=(b_1, v_3b_{-1},b_0)$ satisfy $$U^p+V^p+(unit)\lambda^{(m-1)p}W^p=0$$ which contradicts the minimality of $m$ since we have seen that $m\geq 2$. Thus, either we have constructed a cyclic extension of the required type or one of $I_0,I_1$ and $I_{p-1}$ is non-principal. But then again from the result of Class Field theory as mentioned in the previous section we have a cyclic extension as required.

\section{{Transcendental Computation of the Class Number}}

The Dedekind zeta function for a number field $K$ and its associated Euler product expansion is given by, $$\zeta_K(s)=\sum_I\frac{1}{N(I)^s}=\prod_Q\frac{1}{(1-\frac{1}{N(Q)^s}}\prod$$ where the sum runs over all ideals $I$ of $R$ and the product runs over all prime ideals $Q$ of $R$. The two expressions gives us two ways of computing $\lim_{s\rightarrow1}(s-1)\zeta_K(s)$. The left-hand side is expressed in terms of arithmetic invariants and the right-hand side in terms of invariants of the Galois group. The resulting identity will give a way for computing the Class number $h$ of $K$.

The left-hand limit can be found to be, $$\lim_{s\rightarrow 1}(s-1)\sum_I\frac{1}{N(I)^s}=\lim_{r\rightarrow\infty}\frac{\#\{I\mid N(I)\leq r\}}{r}.$$

The set $\{I\mid N(I)\leq r\}$ can be split according to ideal classes, and we compute for each ideal class $C$, $$z(C)=\lim_{r\rightarrow\infty}\frac{\#\{I\in C\mid N(I)\leq r\}}{r}.$$

Fixing an ideal $I_0\in C$, this latter set is a bijection to $\{aR\subset I_0^{-1}\mid N(a)\leq r.N(I_0)^{-1}\}$, where $N(a)$ denotes the modulus of the norm of $a$.

We then have a natural embedding $K\hookrightarrow K\otimes_Q R$. The image of $J=I_0^{-1}$ is a lattice in $K\otimes_Q R$. Let $\Lambda$ denote the image of $J-\{0\}$ in the quotient $S=(K\otimes_QR)^{\ast}/U$ where $U$ is the image of the group of units in $R$ under the above embedding. There is a natural homomorphism $N:S\rightarrow R^\ast$ which restricts the modulus of the norm on the image of $K$. We then obtain a natural bijection between $\{aR\subset I_0^{-1}\mid N(a)\leq r\}$ and $\{I\in \Lambda \mid N(l)\leq r\}$. Let $\Lambda_r$ denote the image of $(1/r)J-\{0\}$ in $S$, then we have a natural bijection between $\{I\in \Lambda \mid N(l)\leq r^d\}$ and $\{I\in \Lambda_r \mid N(l)\leq 1\}$, where $d$ denotes the degree of $K$ over $Q$.

Let $S_{\leq 1}$ denote the locus of $l\in S$ such that $N(l)\leq 1$. Let $\mu$ denote the Haar measure on $K\otimes_Q R$. This is invariant under the section of $U$ and thus gives a measure also denoted by $\mu$ on $S$. Since $J$ is a lattice in $K\otimes_Q R$ we have, $$\lim_{r\rightarrow\infty}\frac{\#\{I\in \Lambda_r \mid N(l)\leq 1\}}{r^d}=\frac{\mu(S_{\leq 1)}}{\mu(K\otimes_Q R/J)}.$$ Moreover, the denominator can be rewritten as, $$\mu(K\otimes_Q R\j)=N(J)\mu(K\otimes_Q R/R).$$

In particular we see that the limit $z(C)$ is independent of the class $C$. Let $(K\otimes_Q R)_1^\ast$ denote the kernel of the norm map. This is a group and thus we have a Haar measure $\nu$ on it. It can be shown that, $$\mu(S_{\leq 1})=\nu((K\otimes_Q R)_1^\ast/U).$$ From the above calculations we can get, $$\lim_{s\rightarrow 1}(s-1).\zeta_K(s)=h.\frac{\nu((K\otimes_Q R)_1^\ast/U)}{\mu(K\otimes_Q R/R)}.$$

This is called the 'Class number formula' for $K$. It can be noted that the denominator can be computed in a closed form in terms of the discriminant $D$ of the field $K$ and the number of pairs of conjugate complex embeddings $r_2$ of $K$.$$\mu(K\otimes_Q R/R)=\frac{1}{2^{r_2}}.\sqrt{\mid D \mid}$$

However, the numerator is in general more complicated as it involves computing the group of units of $K$. It should be noted that there are other versions of the class number formula too.

To expand the above term we restrict to abelian extensions $K$ of $Q$. The product term on the right can be first grouped according to rational primes $$\prod_Q\frac{1}{(1-\frac{1}{N(Q)^s})}=\prod_q\prod_{Q\mid q}\frac{1}{(1-\frac{1}{N(Q)^s})}.$$

Now for each rational prime $q$ which is unramified in $K$ we have $$\prod_{Q\mid q}\frac{1}{(1-\frac{1}{N(Q)^s})}=\prod_\chi\frac{1}{(1-\frac{\chi(q)}{q^s})}$$ where $\chi$ runs over all characters of the Galois group and $\chi(q)=\chi(Frob_q)$ is the value of $\chi$ on a Frobenius element associated with $q$.

We now define Dirichlet L-series and their Euler product formula as follows, $$L(s,\chi)=\sum_n\frac{\chi(n)}{n^s}=\prod_p\frac{1}{(1-\frac{\chi(p)}{p^s})}$$ where we set $\chi(p)=0$ when $\chi$ is ramified in $p$. We also define the additional factor $$F(s)=\prod_{p~ramified}\frac{1}{(1-\frac{1}{p^{f_p}})^{g_p}}$$ where the product runs over all ramified primes and $f_p$ denotes the residue field extension over $p$ and $g_p$ the number of distinct primes in $K$ lying over $p$. The product expansion of $\zeta_K(s)$ becomes $$\zeta_K(s)=F(s).\prod_\chi L(s,\chi).$$ Thus, the computation of the limit can be reduced to the corresponding computation of the Dirichlet L-series. For the case of unit character we get by comparing with the zeta function, $$\lim_{s\rightarrow 1}(s-1)F(s)L(s,1)=1.$$ So the left-hand limit gives us $$\lim_{s\rightarrow 1}(s-1)\zeta_K(s)=\prod_{\chi\neq 1}L(1,\chi).$$

There is a positive integer $m$ such that $\chi$ is determined on classes modulo $m$ and $\chi$ is $0$ on all primes $p$ dividing it, then $m$ is called the conductor of $\chi$. We rewrite the L-function associated with $\chi$ as follows, $$L(s,\chi)=\sum_{x\in(Z/mZ)^\ast}(\chi(x).\sum_{n\equiv x~(mod~m)}\frac{1}{n^s}).$$ This latter sum can be again rewritten using the identity $$\sum_{i=0}^{m-1}\omega^{xi}=0~or~m$$ according as $x$ is not congruent or congruent to $0$ modulo $m$, where $\omega$ is a primitive m-th root of unity. The second sum then becomes $$\sum_{n\equiv x~(mod~m)}\frac{1}{n^s}=\frac{1}{m}\sum_{n=1}^{\infty}\frac{1}{n^s}\sum_{i=0}^{m-1}\omega^{(x-n)i}.$$

Thus we obtain, $$L(s,\chi)=\frac{1}{m}\sum_{i=0}^{m-1}(\sum_{z\in (Z/mZ)^\ast}\chi(x)\omega^{ix}).\sum_{n=1}^{\infty}\frac{\omega^{-in}}{n^s}.$$ 

The expression, $$\tau_i(\chi)=\sum_{x\in(Z/mZ)^\ast}\chi(x)\omega^{ix}$$ is called the Gaussian sum associated with the integer $i$ and the character $\chi$. If $\chi$ is not a unit character then $\tau_0(\chi)=0$. Moreover, if $i\neq 0$ then we have the identity $$\sum_{n=1}^{\infty}\frac{\omega^{-in}}{n}=-\log(1-\omega^{-i}).$$

Hence, we obtain the formula when $\chi$ is not the unit character $$L(1,\chi)=-\frac{1}{m}\sum_{i=1}^{m-1}\tau_i(\chi).\log(1-\omega^{-i}).$$

\section{{Divisibility of the Class Number by $p$}}
Combining the results that we have proved or stated so far we see that any counter-example of Fermat's Last Theorem for a prime $p\geq 5$ leads to a non-trivial representation $$\rho:Gal(\overline{K}/K)\rightarrow F_p$$ which is unramified everywhere. Kummer called primes which admist such representations irregular. he showed that there are indeed such primes like $37$ and hence this particular attempt to proof the Fermat's Last Theorem fails. We shall now show how to check whether a prime is irregular.

We shall apply the results of the previous Section in the special case where $K$ is the prime cyclotomic field and also to the totally real subfield $L$ described in the previous sections.

First we shall use the divisibility of the class number $h$ of $R$ by the class number $h_{+}$ of $S$ to write $h=h_{+}.h_{-}$ for some integer $h_{-}$. Let $W$ denote the finite cyclic group of roots of unity in $K$. Then we have $U=W.U_{+}$, where $U_{+}$ denotes the group of units in $S$ and so $\#(U/U_{+})=\#(W/\{\pm1\})=p$. We have the natural inclusion $L\otimes_Q R\hookrightarrow K\otimes_Q R$ from which we obtain the isomorphism $$(K\otimes_Q R)^\ast_1/(L\otimes_Q R)_1^\ast=(C_1^\ast/R_1^\ast)^{(p-1)/2}$$ since $(p-1)/2$ is the degree of $L$ over $Q$. From this we can deduce that $$\nu((K\otimes_Q R/R)_1^\ast/U)=\frac{1}{p}.\nu(C_1^\ast/R_1^\ast)^{(p-1)/2}.\nu((L\otimes_Q R)_1^\ast/U_{+}).$$

The formula for computing the discriminant yields $$\mu(K\otimes_Q R/R)=\mu(l\otimes_Q R/S)^2.p^{1/2}$$ since $p$ is the norm of the relative discriminant. Thus, the class number formulas for $K$ and $L$ then give a formula for $h_{-}$, $$\frac{h_{-}.\nu(C_1^\ast/R_1^\ast)^{(p-1)/2}}{p^{3/2},\mu(L\otimes_Q R/S)}=\prod_{\chi(-1)=-1}L(1,\chi).$$

hence, $h_{-}$ can be computed explicitly and in a closed form. in particular the divisibility of $h_{-}$ by $p$ is an easily computable criterion.

The divisibility of $h_{+}$ by $p$ is however more complicated. As remarked earlier the term $\nu((L\otimes_Q R)_1^\ast/U_{+})$ is difficult to compute. However, we have the subgroup $U_{+,cycl}=U_{+}\cap U_{cycl}$ and we can compute $\nu((L\otimes_Q R)_1^\ast/U_{+,cycl})$. In fact we can show that, $$\nu((L\otimes_Q R)_1^\ast/U_{+,cycl})=\mu(L\otimes_Q R/S).\prod_{\chi~even}L(1,\chi)$$ where the product runs over all non-trivial characters $\chi$ such that $\chi(-1)=1$. The class number formula for $h_{+}$ becomes, $$h_{+}=[U_{+}:U_{+,cycl}]=[U:U_{cycl}].$$ This is the first coincidence that makes Kummer's calculations possible.

From the above identity we see that if $p$ divides $h_{+}$ then we have a real unit $u$ such that its p-th power is a cyclotomic unit but $u$ is not itself cyclotomic. Hence $v=u^p$ is a cyclotomic unit which is congruent to an integer modulo $pS$. If we find a $w\in U_{cycl}$ such that $v=w6p$ then we can show easily that $u$ itself is a cyclotomic unit. Let $Q$ denote the quotient group $(S/pS)^\ast(Z/pZ)^\ast$. We obtain a natural homomorphism $$m:U_{cycl}\otimes(Z/pZ)\rightarrow Q$$ which is represented by a square matrix with entries from $F_p$. The preceeding remarks imply that $p\mid h_{+}$ only if $det(m)=0$. The second coincidence that makes Kummer's calculation work is that $det(m)\equiv h_{-}~(mod~p)$.

Thus we see that $p \mid h$ iff $p\mid h_{-}$. Hence we can easily check which primes are regular.

\section{{Acknowledgment}}

I wholeheartedly thank my guide, Prof.~Kapil Hari Paranjape for his guidance and help during the project period. I would also like to thank the Indian Academy of Sciences, Bangalore and the Indian Institute of Science Education and Research (IISER), Mohali for the excellent facilities provided to me during my stay. My thanks are also due to Mr.~Rajib Haloi (IIT, Kanpur) for pointing me towards \cite{prasad} and to my friend, Akshay Kumar Singh (IISER, Kolkata) for providing me with a few of the references below.

\end{document}